\documentclass[fleqn,11pt,a4paper,draft]{article}

\ifx\xyloaded\undefined
  \input{xy}
  \xyoption{all}
  \xyoption{ps}
  \xyoption{dvips}
\else
\fi

\usepackage{fullpage}
\usepackage{style}

\def\first{\hskip1cm&\hskip-1cm}

\DeclareMathOperator{\SL}{SL}
\DeclareMathOperator{\GL}{GL}
\DeclareMathOperator{\Aut}{Aut}
\DeclareMathOperator{\Ext}{Ext}
\DeclareMathOperator{\Tor}{Tor}
\DeclareMathOperator{\gr}{gr}
\DeclareMathOperator{\res}{res}
\DeclareMathOperator{\ind}{ind}
\def\Sp{\mathrm{Sp}} 

\def\HH{H\!H}
\let\fn\rightarrow
\let\ot\otimes
\let\w\omega
\let\l\lambda
\def\C{{\mathbb C}}
\def\N{{\mathbb N}}
\def\Q{{\mathbb Q}}
\def\ZZ{{\mathbb Z}}

\let\eps\varepsilon
\let\T\Tilde
\let\d\partial
\let\*\bullet
\let\epi\twoheadrightarrow

\def\cf{\textit{cf.}}
\def\ie{\textit{i.e.}}
\def\br#1{\langle #1\rangle}

\def\Z{\mathscr{Z}}
\def\B{\mathscr{B}}
\def\G{\mathscr{G}}

\title{Algebra structure on the Hochschild cohomology of the ring of
invariants of a Weyl algebra under a finite group}

\author{Mariano Suarez Alvarez%
\thanks{Departamento de Matem\'atica,
Facultad de Ciencias Exactas, Ingenier\'\i a y Agrimensura, Universidad
Nacional de Rosario. Pellegrini 250. Rosario (2000), Argentina. e-mail:
\texttt{mariano@fceia.unr.edu.ar}.\newline
\hspace*{\parindent} % why on earth isn't \thanks defined \long?
This work was supported by a grant from UBACYT TW69, the
international cooperation project SECyT-ECOS A98E05, and a CONICET 
scholarship.}} 

\date{}

\begin{document}

\maketitle

\begin{abstract}
Let $A_n$ be the $n$-th Weyl algebra, and let
$G\subset\Sp_{2n}(\C)\subset\Aut(A_n)$ be a finite group of linear
automorphisms of $A_n$. In this paper, we compute the multiplicative
structure on the Hochschild cohomology $\HH^\*(A_n^G)$ of the algebra of
invariants of $G$. We prove that, as a graded algebra, $\HH^\*(A_n^G)$ is
isomorphic to the graded algebra associated to the center of the group
algebra $\C\,G$ with respect a filtration defined in terms of the defining
representation of $G$.
\end{abstract}

%%%%%%%%%%%%%%%%%%%%%%%%%%%%%%%%%%%%%%%%%%%%%%%%%%%%%%%%%%%%%%%%%%%%%%
\section{Introduction}\
\label{sect:intro}

\paragraph Let us fix an algebraically closed ground field $\C$ of 
characteristic zero.

\paragraph For $n\in\N$, the $n$-th Weyl algebra $A_n$  is the one freely
generated by elements $p_i$ and $q_i$, $1\leq i\leq n$, subject to the
commutation relations of Heisenberg:
  \begin{alignat*}{2}
  &[p_i,p_j]=[q_i,q_j]=0, &\quad&\forall\,i, j; \\
  &[q_i,p_i]=1,&&\forall\,i;\\
  &[q_j,p_i]=0,&&\text{$\forall\,i, j$ such that $i\neq j$.}
  \end{alignat*}
It can be realized either as the algebra of algebraic differential
operators on the affine space $\mathbb A_n$, or as the Sridharan twisted
enveloping algebra $\mathcal U_f\mathfrak g$ of an abelian Lie algebra
$\mathfrak g$ of dimension $2n$ with respect to any non-degenerate
Chevalley-Eilenberg $2$-cocycle on $\mathfrak g$. It is a simple, left and
right Noetherian algebra of Gabriel-Rentschler Krull dimension $n$,
Gel'fand-Kirillov dimension $2n$ and global homological dimension $n$.

\paragraph Sridharan \cite{Sridharan} shows that the Hochschild cohomology
$\HH^\*(A_n)=H^\*(A_n,A_n)\cong\C$ and,
in fact, that this characterizes the Weyl algebras among the twisted
enveloping algebras of abelian Lie algebras. This result can be interpreted
as the Poincar\'e  lemma for quantum differential forms.
The same methods can be used to show that, dually,
$\HH_\*(A_n)=H_\*(A_n,A_n)\cong\C$,
concentrated in degree $2n$.

\paragraph Consider a finite subgroup $G\subset\Aut(A_1)$ and the
corresponding algebra of invariants $A_1^G$. As $G$ varies, we obtain in
this way a family of algebras, all of which are simple, left and right
Noetherian, with Gel'fand-Kirillov dimension $2$, Krull dimension $1$ and
global homological dimension $1$; in particular, these numeric invariants
do not allow us to separate them. 

Alev and Lambre \cite{AlevLambre} compute the $0$-degree Hochschild
homology of these algebras: they show that $\HH_0(A_1^G)$ is a vector space
of dimension $s(G)-1$, with $s(G)$ the number of irreducible
representations of $G$.  A theorem of Alev \cite{Alev} which describes
$\Aut(A_1)$ implies that each of its finite groups is conjugate to a
subgroup in $\SL_2(\C)\subset\Aut(A_1)$, and the classification up to
conjugation of these is classical; with this information one can compute
$s(G)$ for each of the possible groups, and conclude that the algebras
under consideration are in fact non-isomorphic in pairs, apart from a few
exceptions.

\paragraph If one considers more generally the algebras $A_n^G$ of
invariants of $A_n$ under the action of a finite subgroup
$G\subset\Sp_{2n}(\C)\subset\Aut(A_n)$---we restrict our attention to linear
automorphisms because we have no description of the whole automorphism
group in this case---we again obtain a family of algebras
indistinguishable, for fixed $n$, on the basis of the above numerical
invariants alone. Alev, Farinati, Lambre and Solotar \cite{AFLS} obtain a
generalization of the above formula for $\HH_0(A_1^G)$: they show that
$\nu_k=\dim_\C \HH_k(A_n^G)$ is the number of conjugacy classes of
$G$ whose elements have unity as an eigenvalue with multiplicity exactly
$k$.

We thus see that in general homology is not enough to separate this
algebras, at least without further analysis: one can easily show that these
numbers $\nu_k$ can be computed in terms of the character of the defining
representation and the power maps of the group $G$, and it is known
\cite{Dade} that there are pairs of non-isomorphic finite groups for which
these data coincide.

\paragraph In part, the interest of these computations comes from the wish
to understand the Poisson structures underlying the objects under
consideration.

The algebra $A_n$ has a natural filtration, that of Bernstein, such that
the associated graded object $\gr A_n$ is a polynomial algebra on $2n$
variables, canonically endowed with a Poisson bracket deduced from the
commutator of $A_n$. The action of $\Sp_{2n}(\C)$ on $\gr A_n$ respects
this structure, so that $(\gr A_n)^G$, for $G\subset\Sp_{2n}(\C)$, is
naturally a Poisson algebra. Moreover, one can show that if
$G\subset\Sp_{2n}(\C)$ is a finite subgroup, the graded algebra associated
to $A_n^G$ with respect to the restricted Bernstein filtration is exactly
$(\gr A_n)^G$, with the same Poisson structure. 

In particular, this means that we can regard the algebras $A_n^G$ as
quantizations of the Poisson algebras $(\gr A_n)^G$. This is the point of
view of \cite{AlevLambre2}, where the authors show that, in a precise
sense, the $0$-degree Hochschild homology of $A_1^G$ approximates the
$0$-degree Poisson homology of $(\gr A_1)^G$. This idea cannot be
transferred to the general case because we have no good definition of
Poisson homology for non-smooth algebras and $(\gr A_n)^G$ is non-smooth
for most finite subgroups $\Sp_{2n}(\C)$.

\paragraph It is a result of \cite{AFLS} that there is a duality between
the homology and the cohomology of the algebras at hand. In particular, we
have $\dim_\C \HH^k(A_n^G)=\nu_{2n-k}$ for each finite subgroup
$G\subset\Sp_{2n}(\C)$. In this paper we complete the computation of
Hochschild cohomology making the algebra structure on $\HH^\*(A_n^G)$
explicit. The final result is the following:

\begin{Theorem}\label{the:theorem}
Let $G\subset\Sp_{2n}(\C)$ be a finite subgroup. Let $G$ act naturally on
the $n$-th Weyl algebra $A_n$. The subspace $V\subset A_n$ 
spanned by the standard generators is $G$-invariant for this action. Define
$d:G\fn\N_0$ by $d(g)=2n-\dim_\C V^g$. For each $p\geq0$, write $F_p\C\,G$
the subspace of the group algebra $\C\,G$ spanned by the elements $g\in G$
such that $d(g)\leq p$.

$F_\*\C\,G$ is an algebra filtration on $\C\,G$, so it restricts to an
algebra filtration on the center $\Z\,G$ of $\C\,G$. There is an
algebra isomorphism $\HH^\*(A_n^G)\cong\gr\Z\, G$.
\end{Theorem}

\paragraph It is very easy to construct examples of the situation
considered in the theorem. If $G$ is a finite group and $V$ is a faithful
$G$-module of degree $n$, $G$ acts faithfully on the algebra of algebraic
differential operators on $V$, which is isomorphic to $A_n$, so we can
regard $G\subset\Aut(A_n)$. One sees that---using the notation of the
theorem---$d(g)$ is simply two times the codimension of the subspace of $V$
fixed by $g$.

One particularly nice instance of this arises when we consider the
canonical action of the Weyl group corresponding to a Cartan subalgebra
$\mathfrak h$ of a semi-simple Lie algebra $\mathfrak g$ on the algebra of 
regular differential operators on $\mathfrak h^*$, the dual space of
$\mathfrak h$.

\paragraph In the next section, we recall the construction of the
multiplicative structure on the Hoch\-schild theory, and indicate the
reductions leading to its determination in our particular case. In section
\ref{sect:computation} we carry out the various explicit computations
needed for the proof of the theorem.

%%%%%%%%%%%%%%%%%%%%%%%%%%%%%%%%%%%%%%%%%%%%%%%%%%%%%%%%%%%%%%%%%%%%%%
\section{The multiplicative structure on $\HH^\*(A_n^G)$}
\label{sect:mult}

\paragraph Let us fix from now on $n\in\N$ and a finite subgroup
$G\subset\Sp_{2n}(\C)$. We consider the natural action of $G$ on the $n$-th
Weyl algebra $A_n$ by linear automorphisms.

\paragraph The computation of $\HH^\*(A_n^G)$ presented in \cite{AFLS} is
based on the fact that $A_n^G$ and the crossed product $A_n\rtimes G$ are
Morita equivalent; since Hochschild cohomology groups are invariant under
this kind of equivalence, in order to determine $\HH^\*(A_n^G)$, one can instead
choose to compute $\HH^\*(A_n\rtimes G)$. Now, the algebra
structure on $\HH^\*(A_n^G)$ can be defined in terms of the composition of
iterated self-extensions of $A_n^G$ in the category of $A_n^G$-bimodules;
since this procedure is clearly invariant under equivalences, the
\emph{algebras} $\HH^\*(A_n^G)$ and $\HH^\*(A_n\rtimes G)$ coincide.

\paragraph The next reduction depends on results of Stefan \cite{Stefan}
and others, which show that, in our situation, for each $A_n\rtimes G$-bimodule
$M$ there is a natural spectral sequence with initial term $E_2^{p,q}\cong
H^p(G,H^q(A_n,M))$ converging to $H^\*(A_n\rtimes G,M)$. Since our ground
field has characteristic zero, group cohomology is trivial in positive
degrees, and this spectral sequence immediately degenerates, giving us
natural isomorphisms $H^\*(A_n\rtimes G,M)\cong H^\*(A_n,M)^G$.

We set $M=A_n\rtimes G$. It is easy to see that Stefan's spectral sequence
is a spectral sequence of algebras in this case---for example, by using the
resolutions given by the bar construction in order to compute the
cohomologies of $G$ and of $A_n$. The distribution of zeros in its initial
term implies that there are no extension problems, neither in order to
compute the cohomology groups nor to compute the product maps. We thus
conclude that the isomorphism between $\HH^\*(A_n\rtimes G)$ and
$H^\*(A_n,A_n\rtimes G)^G$ is multiplicative; we will determine this last
algebra.

\paragraph Let us recall the construction of the multiplicative structure
on the functor $H^\*(A_n,-)=\Ext_{A_n^e}^\*(A_n,-)$. We fix a projective
resolution $X^\*\epi A_n$ of $A_n$ as a $A_n^e$-module. Since plainly
$\Tor^{A_n}_p(A_n,A_n)=0$ for $p\geq0$, $X^\*\otimes_{A_n}X^\*$ is an
acyclic complex over $A_n\ot_{A_n}A_n\cong A_n$. It follows from \cite{CE},
proposition  \textsc{IX.2.6}, that it is a projective resolution of $A_n$ as an
$A_n^e$-module.

In particular, there is a morphism $\Delta:X^\*\fn X^\*\ot_{A_n}X^\*$ of
resolutions lifting the identity map of $A_n$. If $M$ and $N$ are
$A_n^e$-modules, the product
  \begin{equation}\label{cup0}
  \cup:\Ext_{A_n^e}^\*(A_n,M)\ot\Ext_{A_n^e}^\*(A_n,N)\fn
  	\Ext_{A_n^e}^\*(A_n,M\ot_{A_n}N)
  \end{equation}
is induced by the composition
  \[\xymatrix@R-10pt{
  \hom_{A_n^e}(X^\*,M)\ot\hom_{A_n^e}(X^\*,N)\ar[d]^-\psi \\
  \hom_{A_n^e}(X^\*\ot_{A_n}X^\*,M\ot_{A_n}N)\ar[d]^-\Delta \\
  \hom_{A_n^e}(X^\*,M\ot_{A_n}N)
  }\]
with $\psi$ standing for the evident $\hom$-$\ot$ ``interchange map'', up to the
canonical isomorphisms
  \[
  H(\hom_{A_n^e}(X^\*,M))\ot H(\hom_{A_n^e}(X^\*,N))\cong
  H(\hom_{A_n^e}(X^\*,M)\ot \hom_{A_n^e}(X^\*,N)).
  \]

\paragraph When $M=N=A_n\rtimes G$, we can compose the map \eqref{cup0}
with the morphism induced by the product $\mu:(A_n\rtimes
G)\ot_{A_n}(A_n\rtimes G)\fn A_n\rtimes G$. We obtain in this way the
internal product of $H^\*(A_n,A_n\rtimes G)$.

The additivity of the functors involved and the decomposition $A_n\rtimes
G\cong\bigoplus_{g\in G}A_ng$ of $A_n\rtimes G$ as an $A_n^e$-module---here
and elsewhere $A_ng$ is the $A_n^e$-module obtained from $A_n$ by twisting
the right action by the automorphism $g$---have the consequence
that this product is determined by its restrictions
  \[
  \cup:\Ext_{A_n^e}^\*(A_n,A_ng)\ot\Ext_{A_n^e}^\*(A_n,A_nh)\fn
  	\Ext_{A_n^e}^\*(A_n,A_ngh),
  \]
which we will compute in the next section.

%%%%%%%%%%%%%%%%%%%%%%%%%%%%%%%%%%%%%%%%%%%%%%%%%%%%%%%%%%%%%%%%%%%%%%
\section{Explicit computations}
\label{sect:computation}

\paragraph First of all, let us consider a filtered $\C$-algebra $A$ with a
positive ascending filtration such that the associated graded algebra $\gr
A$ has no zero divisors. For each $x\in A$, we denote $s(x)$ the principal
symbol of $x$ in $\gr A$. 

Consider commuting  elements $x_1,\dots,x_n\in A$, and write, for each $k$
such that $0\leq k\leq n$, $I_k$ and $I_k'$ the left ideals generated by
$x_1,\dots,x_k$ and $s(x_1),\dots,s(x_k)$ in $A$ and $\gr A$, respectively;
in particular, $I_0=I_0'=0$. Let us suppose that we have, for $1\leq k\leq
n$, 
  \begin{equation}\label{i}
  a\in\gr A, as(x_k)\in I_{k-1}'\Rightarrow a\in I_{k-1}'.
  \end{equation}
Let $k$ be such that $1\leq k\leq n$, $a\in F_mA$ and suppose $ax_k\in
I_{k-1}$, so that there are $a_i\in A$ for $i=1,\dots,k-1$ with
$ax_k=\sum_{i=1}^{k-1}a_ix_i$. Then
$s(ax_k)=s(\sum_{i=1}^{k-1}a_ix_i)\in I_{k-1}'$, and, as 
$s(ax_k)=s(a)s(x_k)$ because $\gr A$ is a domain, we see from \eqref{i}
that $s(a)\in I_{k-1}'$. We conclude that there exist $b_i\in A$, for each 
$i=1,\dots,k-1$, and $a'\in F_{m-1}A$ such that $a=\sum_{i=1}^{k-1}b_ix_i+a'$. 
We have
  \[
  a'x_k=ax_k-\sum_{i=1}^{k-1}b_ix_ix_k=ax_k-\sum_{i=1}^{k-1}b_ix_kx_i\in
  I_{k-1};
  \]
by induction, this tells us that $a'\in I_{k-1}$, and, consequently, that
$a\in I_{k-1}$. 

We have shown that
  \begin{equation}\label{ii}
  a\in A, ax_k\in I_{k-1}\Rightarrow a\in I_{k-1},
  \end{equation}
for each $k$ such that $1\leq k\leq n$, and we see that this is a condition
that can be tested up to an appropriate filtration.

\paragraph Let $A_n$ be the $n$-th Weyl algebra over $\C$, with generators
$p_i$, $q_i$, for $i=1,\dots,n$. Let $\mu:A_n^e\fn A_n$ be the canonical
augmentation, and put $I=\ker\mu$. If $x\in A$, we'll write $\d x=1\ot
x-x\ot 1$; plainly, $\d x\in I$. A trivial computation shows that if $x,y\in
A_n$ are such that $[x,y]\in \C1$, $[\d x,\d y]=0$. In particular, 
the elements $\d p_i$ and $\d q_i$, $i=1,\dots,n$, which generate $I$ as a
left $A_n^e$-module, commute.

There is an isomorphism of algebras $\phi:A_n^e\fn A_{2n}$ uniquely
determined by the conditions
  \[
  \phi(p_i\ot1)=p_i,\qquad
  \phi(q_i\ot1)=q_i,\qquad
  \phi(1\ot p_i)=q_{i+n},\qquad\text{and}\qquad
  \phi(1\ot q_i)=p_{i+n},
  \]
for each $i=1,\dots n$. One has $\phi(\d p_i)=-p_i+q_{i+n}$ and
$\phi(\d q_i)=-q_i+p_{i+n}$, and these elements obviously commute.
When we consider the Bernstein filtration on $A_{2n}$, $\gr A_{2n}$ turns
out to be a polynomial algebra on variables $x_i=s(p_i)$ and
$y_i=s(q_i)$ for $i=1,\dots,2n$; moreover, if $i=1,\dots,n$, 
  \[
  s(\phi(\d p_i))=-x_i+y_{i+n},\qquad 
  s(\phi(\d q_i))=-y_i+x_{i+n}.
  \]
It is clear that in this case \eqref{i} is satisfied, so that in turn the elements
$\d p_1,\dots,\d p_n$ and $\d q_1,\dots,\d q_n$ satisfy \eqref{ii}.

We see that the left augmented algebra $\mu:A_n^e\fn A_n$ satisfies the
hypotheses of theorem \textsc{VIII.4.2} of \cite{CE}; in particular, if we
let $V=\bigoplus_{i=1}^n\C\,p_i\oplus\bigoplus_{i=1}^n\C\,q_i$, we have a
projective resolution of $A_n$ as a left $A_n^e$-module of the form 
$A_n^e\ot \Lambda^\bullet V\epi A_n$ with differentials
$d:A_n^e\ot\Lambda^pV\fn A_n^e\ot\Lambda^{p-1}V$ given by
  \[
  d( a\ot v_1\wedge\cdots\wedge v_p)
    = \sum_{i=1}^p(-1)^{i+1} a\d v_i\ot v_1\wedge\cdots\wedge\Hat
    v_i\wedge\cdots\wedge v_p.
  \]
There are, accordingly, natural isomorphisms
  \[
  \Ext_{A_n^e}^\*(A_n,M)
    \cong H(\hom_{A_n^e}(A_n^e\ot\Lambda^\* V,M))
    \cong H(\hom(\Lambda^\* V,M))
  \]
for each left $A_n^e$-module $M$, where, in the last term, homology is
computed with respect to differentials
$d:\hom(\Lambda^dV,M)\fn\hom(\Lambda^{d+1}V,M)$ such that
  \begin{align*}\first
  df(v_1\wedge\cdots\wedge v_{p+1}) =
    \sum_{i=1}^{p+1}(-1)^{i+1}\partial v_if(v_1\wedge\cdots\wedge\Hat
    v_i\wedge\cdots\wedge v_{p+1}) \\
    &= \sum_{i=1}^{p+1}(-1)^{i+1}[v_i,f(v_1\wedge\cdots\wedge\Hat
    v_i\wedge\cdots\wedge v_{p+1})]
  \end{align*}
for each $f:\Lambda^d V\fn M$.

\paragraph\label{constr} Keeping the notations introduced in the previous
paragraph, let $g\in\Sp(V)$, where we view $V$ as a symplectic space in the
usual way; there is a unique decomposition $V=V_1^g\oplus V_2^g$ preserved
by $g$ and such that $g|V_1^g=Id_{V_1^g}$ and
$Id_{V_2^g}-g|_{V_2^g}\in\GL(V_2^g)$.   Let $d(g)=\dim V_2^g$.
Whenever possible, we will suppress reference to the automorphism $g$ in
our notation.

Let $\w\in(\Lambda^dV_2)^*\setminus 0$. The decomposition
$V=V_1\oplus V_2$ induces a decomposition
  \[
  \Lambda^dV=\bigoplus_{p+q=d}\Lambda^pV_1\otimes\Lambda^qV_2;
  \]
in particular, we see that $\Lambda^dV_2$ can be identified with a subspace
of $\Lambda^dV$, and, in this identification, admits a natural complement; 
we extend $\w$ to the whole of $\Lambda^dV$ prescribing it to be zero on
this complement.

We define $\T\w:\Lambda^dV\fn A_ng$ by setting $\T\w(v)=\w(v)g$. 
We will show that $\T\w$ represents a non-zero homology class of degree $d$
in the complex $\hom(\Lambda^\bullet V,Ag)$ considered above.

Let us fix a basis $v_1,\dots,v_{2n}$ of $V$ an such a way that $v_1,\dots,v_d\in
V_2$ and $v_{d+1},\dots,v_{2n}\in V_1$, and choose indices $1\leq
i_1<\dots<i_{d+1}\leq 2n$; we have
  \begin{align*}\first
  d\T\w(v_{i_1}\wedge\cdots\wedge v_{i_{d+1}})
    =
    \sum_{j=1}^{d+1}(-1)^{j+1}[v_{i_j},\T\w(v_{i_1}\wedge\cdots\wedge
    \Hat v_{i_j}\wedge\cdots\wedge v_{i_{d+1}})] \\
    &= 
    \sum_{j=1}^{d+1}(-1)^{j+1}[v_{i_j},\w(v_{i_1}\wedge\cdots\wedge
    \Hat v_{i_j}\wedge\cdots\wedge v_{i_{d+1}})]_gg;
  \end{align*}
in the second equality, we are seeing $\w$ as taking values in $A_n$, and
writing $[x,y]_g=xy-yg(x)$.

When $i_d>d$, $\w(v_{i_1}\wedge\cdots\wedge 
\Hat v_{i_j}\wedge\cdots\wedge v_{i_{d+1}})=0$ for every $j$ such that
$1\leq j\leq d+1$, and this implies that, in this case, 
$d\T\w(v_{i_1}\wedge\cdots\wedge v_{i_{d+1}})=0$; if, on the contrary,
$i_d\leq d$, necessarily $i_j=j$ for each $1\leq j\leq d$, and
$\w(v_{i_1}\wedge\cdots\wedge 
\Hat v_{i_j}\wedge\cdots\wedge v_{i_{d+1}})=0$ if $1\leq j\leq d$, so 
we have simply
  \[
  d\T\w(v_{i_1}\wedge\cdots\wedge v_{i_{d+1}})
    = (-1)^d[v_{i_{d+1}},\w(v_1\wedge\cdots\wedge v_d)]_gg.
  \]
As the values of $\w$ are in the center of $A_n$ and $v_{i_{d+1}}\in V_1$,
this twisted commutator vanishes, and, again, we have
$d\T\w(v_{i_1}\wedge\cdots\wedge v_{i_{d+1}})=0$. Having considered each
element in a basis of $\Lambda^{d+1}V$, we conclude $d\T\w=0$, i.e.~$\T\w$ 
is a $d$-cocycle.

Suppose now that there is an $h\in\hom(\Lambda^{d-1}V,A)$ such that $d(hg)=\T\w$.
Then, writing $h_i=h(v_1\wedge\cdots\wedge\Hat v_i\wedge\cdots\wedge
v_d)\in A_n$, we have
  \begin{equation}\label{abs}
  \sum_{i=1}^d(-1)^{i+1}[v_i,h_i]_g = \omega(v_1\wedge\cdots\wedge v_d) \in
  \C1.
  \end{equation}
In particular, $[V_2,A_n]_g\cap \C1\neq0$.

For each $i=1,\dots,n$, let $V^i$ be the subspace of $V$ spanned by $p_i$
and $q_i$, and let $A^i$ be the subalgebra of $A_n$ generated by $p_i$ and $q_i$.
Abusing a little of our notation, we see that $V=\oplus_{i=1}^nV^i$ and
$A_n=\otimes_{i=1}^nA^i$.  Suppose, without any loss of generality, that
$V_2$ and $V_1$ are generated by $p_i$, $q_i$ for $i=1,\dots,d$ and for
$i=d+1,\dots,n$ respectively, and that each $V^i$ is preserved by $g$; let
us write $g_i=g|_{V^i}$. We have
  \begin{equation}\label{iii}
  [V_2,A_n]_g = \sum_{i=1}^d[V^i,A_n]_g
    = \sum_{i=1}^d A^1\ot\cdots\ot[V^i,A^i]_{g_i}\ot\cdots\ot A^n.
  \end{equation}
Theorem 4 in \cite{AlevLambre} states, among other things, that, for $g$ an
automorphism of $A_1$ different from the identity,
$A_1=\C1\oplus[A_1,A_1]_g$.   This implies that
$A^i=\C1\oplus[A^i,A^i]_{g_i}$, and, using \eqref{iii}, that $\C1$ and
$[V_2,A_n]_g$ are transversal subspaces in $A_n$. This contradicts
\eqref{abs}; we have thus proved that $\T\w$ cannot be a coboundary in our
complex.

\paragraph In fact, this construction produces \emph{all} $d$-cocycles with
values in $\C g\subset A_ng$. To see this, let $\eta:\Lambda^dV\fn A_ng$ be
one such cocycle, let $v_1,\dots,v_{2n}$ be a basis of $V$ as above with
respect to which  $g$ acts diagonally; for each $1\leq i\leq 2n$, 
let $\lambda_i\in\C$ be such that $gv_i=\lambda_iv_i$. If
$1\leq i_1<\dots<i_{d+1}\leq 2n$ and $s$ satisfies $i_s\leq d<i_{s+1}$, we
have that
  \begin{align*}
  0 &= d\eta(v_{i_1}\wedge\dots\wedge v_{i_{d+1}}) \\
    &= \sum_{j=1}^{d+1}(-1)^{j+1} [v_{i_j},
	 \eta(v_{i_1}\wedge\dots\wedge\Hat v_{i_j}\wedge\dots\wedge v_{i_{d+1}})]_g \\
    &= \sum_{j=1}^{d+1}(-1)^{j+1} [v_{i_j},
	 \eta(v_{i_1}\wedge\dots\wedge\Hat v_{i_j}\wedge\dots\wedge
	 v_{i_{d+1}})] \\
    &
       + \sum_{j=1}^s(-1)^{j+1}(1-\lambda_{i_j})
	 \eta(v_{i_1}\wedge\dots\wedge\Hat v_{i_j}\wedge\dots\wedge
	 v_{i_{d+1}})v_{i_j};
  \end{align*}
since the first sum in the last member is zero, we see that, for $1\leq
j\leq s$, 
  \[
  \eta(v_{i_1}\wedge\dots\wedge\Hat v_{i_j}\wedge\dots\wedge v_{i_{d+1}})=0.
  \]

If now $1\leq r_1<\dots<r_d\leq 2n$ and $r_d>d$, there exists
$t\in\{1,\dots,d\}\setminus\{r_1,\dots,r_d\}$; let $1\leq
i_1<\dots<i_{d+1}\leq 2n$ be such that
$\{i_1,\dots,i_{d+1}\}=\{r_1,\dots,r_d\}\cup\{t\}$ and suppose
$i_f=t$. Our previous observation implies that
  \[
  \eta(v_{r_1}\wedge\dots\wedge v_{r_d})
  = \eta(v_{i_1}\wedge\dots\wedge\Hat v_{i_f}\wedge\dots\wedge v_{i_{d+1}})
  = 0.
  \]
We see that $\eta(v_{r_1}\wedge\dots\wedge v_{i_d})$ vanishes unless
$r_d\leq d$, and in this case $v_{r_j}=v_j$ if $1\leq j\leq d$. It is clear
now that $\eta$ is one of the cocycles constructed the previous
paragraph.

\paragraph In our situation, and using the notations of the end of
\pref{constr}, we have an iterated product, c.f.~\cite{CE}~\textsc{XI.1},
  \[
  \vee:\bigotimes_{i=1}^n\Ext_{(A^i)^e}^\*(A^i,A^ig_i)
  \fn \Ext_{A_n^e}^\*(A_n,A_ng)
  \]
which is an isomorphism in view of theorem \textsc{XI.3.1}, \textit{loc.cit.} 
On the other hand, we know from \cite{AFLS} that, for an algebra automorphism $g$ of $A_1$,
  \[
  \dim_\C\Ext_{A_1^e}^\*(A_1,A_1g)=\left\{
    \begin{array}{ll}
    1, & \text{if $g=Id_{A_1}$ and $\bullet=0$ or if $g\neq Id_{A_1}$ and $\bullet=2$;}\\
    0, & \text{in any other case.}
    \end{array}
    \right.
  \]
From these two facts, we easily deduce that $\Ext_{A_n^e}^\bullet(A_n,A_ng)$
is trivial except in degree $d$, where it is one-dimensional.
Comparing dimensions, we see that the map $\w\mapsto\T\w$ is an 
isomorphism
  \begin{equation}\label{iso:1}
  \Lambda^d(V_2)^*[-d]\cong \Ext_{A_n^e}^\bullet(A_n,A_ng).
  \end{equation}
Here $M[-d]$ denotes the $d$-th suspension of a graded space $M$.

\paragraph\label{G} Let $G\subset\Sp(V)$ be a finite subgroup. The natural
action of $G$ on $V$ extends to an homogeneous action on the exterior
algebra $\Lambda^\bullet V$, on one side, and, on the other, to an action
by algebra automorphisms on $A_n$, and, thence, on $A_n^e$. With respect
to these actions, each module in the resolution $A_n^e\ot\Lambda^\bullet
V\epi A_n$ is a $G$-module, and the differentials are $G$-equivariant.

If $M$ is a left $A_n^e\rtimes G$-module, there is an homogeneous action
of $G$ on $\hom_{A_n^e}(A^e_n\ot\Lambda^\bullet V,M)$, which is natural with
respect to morphisms $M\fn M'$ of $A_n^e\rtimes G$-modules, and, under the
isomorphism of $\C$-spaces 
$\hom_{A_n^e}(A^e_n\ot\Lambda^\bullet V,M)\cong\hom(\Lambda^\bullet V,M)$,
corresponds to the usual diagonal action of $G$.
Passing to homology, we obtain an action of $G$ on 
$\Ext_{A_n^e}^\bullet(A_n,M)\cong H(\hom(\Lambda^\bullet V,M))$.

In particular, $\Ext_{A_n^e}^\bullet(A_n,A_n\rtimes G)$ is, in a natural
way, a graded $G$-module. In view of the decomposition
$A_n\rtimes G\cong\oplus_{g\in G}A_ng$ of $A_n\rtimes G$ as a left $A_n^e$-module 
and the considerations of the previous paragraph, we have an isomorphism
  \[
  \Ext_{A_n^e}^\bullet(A_n,A_n\rtimes G)\cong
  \bigoplus_{g\in G}(\Lambda^{d(g)}V_2^g)^*g[-d(g)],
  \]
With respect to this isomorphism, the action of $G$ can be described in the
following way: let $g, h\in G$ and 
$\w\in(\Lambda^{d(g)}V_2^g)^*$; left multiplication by $h^{-1}$ induces an
isomorphism $V_2^{hgh^{-1}}\fn V_2^g$ which, in turn, determines an
isomorphism
$h^\flat:(\Lambda^{d(g)}V_2^g)^*\fn(\Lambda^{d(hgh^{-1})}V_2^{hgh^{-1}})^*$;
in this notation, we have
  \[
  h(\w g)=h^\flat(\w) hgh^{-1}.
  \]
The verification of this claim reduces to a simple computation.

\paragraph Since obviously $\Tor^{A_n}_+(A_n,A_n)=0$, we know that
  \[
  A_n\ot\Lambda^\bullet V\ot A_n\ot\Lambda^\bullet V\ot A_n\cong
  (A_n^e\ot\Lambda^\bullet V)\ot_{A_n} (A_n^e\ot\Lambda^\bullet V)
  \epi A_n\ot_{A_n} A_n \cong A_n
  \]
is a projective resolution of $A_n$ as a left $A_n^e$-module. There is a
morphism of resolutions of $A_n$ over $1_{A_n}$,
$\Delta:A_n\ot\Lambda^\bullet V\ot
A_n\fn A_n\ot\Lambda^\bullet V\ot A_n\ot\Lambda^\bullet V\ot A_n$, given, in
each degree $d$, by
  \[
  \Delta(a\ot v_1\wedge\cdots\wedge v_d\ot b)
  = \sum_{\substack{p+q=n\\(i,j)\in S_{p,q}}}
    \eps(i,j)\,(a\ot v_{i_1}\wedge\cdots\wedge v_{i_p}\ot 
             1\ot v_{j_1}\wedge\cdots\wedge v_{j_q}\ot b),            
  \]
if we let $S_{p,q}$ be the set of $(p,q)$-shuffles in the symmetric group
$S_{p+q}$, and, for each such shuffle $(i,j)$, $\eps(i,j)$ is the signature of
$(i,j)$.

Given $A_n^e$-modules $M$ and $N$, the product
  \begin{equation}\label{cup}
  \cup:\Ext_{A_n^e}^\bullet(A_n,M)\ot\Ext_{A_n^e}^\bullet(A_n,N)
       \fn \Ext_{A_n^e}^\bullet(A_n,M\ot_{A_n} N)
  \end{equation}
is as explained in section \ref{sect:mult}.
Explicitly, under the usual identifications,
if $\xi\in\hom(\Lambda^p V,M)$ and $\zeta\in\hom(\Lambda^q V,N)$, the
product
$\xi\cup\zeta\in\hom(\Lambda^{p+q}V,M\ot_{A_n} N)$ is such that
  \[
  (\xi\cup\zeta)(v_1\wedge\cdots\wedge v_{p+q})
  = \sum_{(i,j)\in S_{p,q}}
    \eps(i,j)\;\xi(v_{i_1}\wedge\cdots\wedge v_{i_p})\ot
             \zeta(v_{j_1}\wedge\cdots\wedge v_{j_q}).
  \]
If there is a group $G$ acting like in \pref{G}, from general principles
or simply in view of this formula, we know that if $M$ and $N$ are
$A_n^e\rtimes G$-modules, the product \eqref{cup} is  $G$-equivariant.

\paragraph In the situation of \pref{G}, choose an arbitrary
$G$-invariant inner product on $V$. It is easy to see that, for each $g\in
G$, $V_1^g$ and $V_2^g$ are mutual orthogonal complements in $V$. If
$g,h\in G$, we have
  \[
  (V_2^g+V_2^h)^\perp
  = {V_2^g}^\perp\cap{V_2^h}^\perp
  = V_1^g\cap V_1^h
  \subset V_1^{gh}
  \]
so that
  \begin{equation}\label{incl}
  V_2^{gh}={V_1^{gh}}^\perp
  \subset(V_2^g+V_2^h)^{\perp\perp}
  = V_2^g+V_2^h.
  \end{equation}

There is an isomorphism $A_ng\ot_{A_n}A_nh\cong A_ngh$ of 
$A_n^e$-modules, which we regard as an identification. Setting $M=A_ng$ and
$N=A_nh$ in  \eqref{cup}, we have a product map
  \begin{equation}\label{cup:gh}
  \cup:\Ext_{A_n^e}^\bullet(A_n,A_ng)\ot\Ext_{A_n^e}^\bullet(A_n,A_nh)
       \fn \Ext_{A_n^e}^\bullet(A_n,A_ngh)
  \end{equation}
From degree considerations, we see that this is trivial unless
$d(gh)=d(g)+d(h)$; if this is the case, \eqref{incl} implies that
$V_2^{gh}=V_2^g\oplus V_2^h$. Let $\w\in(\Lambda^{d(g)}V_2^g)^*$ and
$\phi\in(\Lambda^{d(h)}V_2^h)^*$ be non-zero forms, and consider a basis
$v_1,\dots,v_{2n}$ of $V$ such that $v_1,\dots,v_{d(g)}$ is a basis of
$V_2^g$, $v_{d(g)+1},\dots,v_{d(g)+d(h)}$ is a basis of $V_2^h$, and
$v_{d(g)+d(h)+1},\dots,v_{2n}$ is a basis of $V_1^{gh}$. Let $1\leq
r_1<\dots<r_{d(g)+d(h)}\leq 2n$ be arbitrary indices. If $r_{d(g)}>d(g)$ then
  \begin{equation}\label{it}
  \begin{split}
  (\T\w\cup\T\phi)(v_{r_1}\wedge&\dots\wedge v_{r_{d(g)+d(h)}})
  =\\
  &\sum_{(i,j)\in S_{d(g),d(h)}}\eps(i,j)
     \T\w(v_{r_{i_1}}\wedge\dots\wedge v_{r_{i_{d(g)}}})
     \T\phi(v_{r_{i_{d(g)}+1}}\wedge\dots\wedge v_{r_{i_{d(g)+d(h)}}})
  \end{split}
  \end{equation}
is zero because, for each $(i,j)\in S_{d(g),d(h)}$, $v_{r_{i_{d(g)}}}>d$, so
the second factor in each term of the sum vanishes. If $r_{d(g)}\leq d$ but
$r_{d(g)+d(h)}>d(g)+d(h)$, a similar reasoning shows \eqref{it} is also zero.

We thus see that unless $r_i=i$ for each $1\leq i\leq d(g)+d(h)$,
$(\T\w\cup\T\phi)(v_{r_1}\wedge\dots\wedge v_{r_{d(g)+d(h)}})=0$; that is,
$\T\w\cup\T\phi$ is one of the cocycles constructed in \pref{constr}. It is
not cohomologous to zero, because it is not zero on
$\Lambda^{d(gh)}V_2^{gh}$, since
  \[
  (\T\w\cup\T\phi)(v_1\wedge\dots\wedge v_{d(g)+d(h)})  
  = \w(v_1\wedge\dots\wedge v_{d(g)})\phi(v_{d(g)+1}\wedge\dots\wedge v_{d(g)+d(h)})
  \neq 0.
  \]

We conclude that \eqref{cup:gh} is either an isomorphism or zero, depending  on
whether $d(gh)=d(g)+d(h)$ or not.

\paragraph We can choose a non zero element
$\w_g\in(\Lambda^{d(g)}V_2^g)^*$ for each $g\in G$ in the following way:
let $\nu\in(\Lambda^2V)^*$ be the symplectic form on $V$; since the action
of $G$ preserves $\nu$, $\nu|_{\Lambda^2V_2^g}$ is a symplectic form on
$V_2^g$ for each $g$; in particular, the $d(g)/2$-th exterior power
$\omega_g=(\nu|_{\Lambda^2V_2^g})^{d(g)/2}\in(\Lambda^{d(g)}V_2^g)^*\setminus0$.
It is clear that when $g, h\in G$ are such that $d(g)+d(h)=d(gh)$,
$\T\w_g\cup\T\w_h=\T\w_{gh}$ because $V_2^g\oplus V_2^h=V_2^{gh}$.
Moreover, these elements are compatible with the action of $G$ on
$\Ext_{A_n^e}^\bullet(A_n,A_n\rtimes G)$, in the sense that
$g\T\w_h=\T\w_{ghg^{-1}}$, because the action is symplectic.

We thus see that in terms of the basis $\{\w_g\}_{g\in G}$ both the
structure constants and the action of $G$ become particularly pleasant.

\paragraph Consider the filtration $F_\*\C\,G$ on $\C G$ such that
$F_p\C\,G$ is spanned by the elements $g\in G$ such that $d(g)\leq p$. 
Equation \eqref{incl} implies that this is an algebra filtration on $\C\,G$.

It is clear from the previous paragraph that the map
  \[
  \T\w_g\in\Ext_{A_n^e}^\*(A_n,A_n\rtimes G)\mapsto s(g)\in\gr\C\, G
  \]
is an algebra isomorphism. The compatibility of the chosen basis of the
domain of this map with the action of $G$ tells us that this map in in fact
$G$-equivariant, and, hence, that there is an isomorphism of graded
algebras $\Ext_{A_n^e}^\*(A_n,A_n\rtimes G)^G\cong(\gr\C\,G)^G$.

\paragraph Write $\Z\,G$ the center of $\C\,G$, and consider on it the filtration
induced by $F_\*\C\,G$. It is clear that $(\C\,G)^G=\Z\,G$, so that 
$\gr\Z\,G\subset(\gr\C\,G)^G$, and in fact this is an equality, because
passing to the associated graded objects preserves the dimension.
This proves theorem \pref{the:theorem}.

%%%%%%%%%%%%%%%%%%%%%%%%%%%%%%%%%%%%%%%%%%%%%%%%%%%%%%%%%%%%%%%%%%%%%% 
\section{Some examples}
\label{sect:examples}

\paragraph As mentioned in the introduction, it is very easy to construct
examples of the situation considered in our theorem \pref{the:theorem}.
Indeed, let $G$ be a finite group and choose a faithful $G$-module $V$ of
degree $n$; $G$ acts faithfully on the algebra of regular algebraic
differential operators on $V$, which is isomorphic to $A_n$, so we can
regard $G\subset\Aut(A_n)$. It is very easy to show that---in the notation
of the theorem---$d(g)$ is simply two times the codimension of the subspace
of $V$ fixed by $g\in G$.

\paragraph Let us write $C[G]$ the algebra of $\C$-valued central
functions on $G$. 
It is well-known---\cf~\cite{Knutson}---that $C[G]$ is canonically
endowed with the structure of a $\lambda$-ring with respect to which the
Adams operations are given by $\psi^k(f)(g)=f(g^k)$ for $k\geq0$,
$f\in C[G]$ and $g\in G$.

Let $t$ be a variable, and let $p\in C[G][t]$ be the central function
with polynomial values such that, for each $g\in G$, $p(g)$ is the
characteristic polynomial of $g$ in the representation $V$; define now
$q\in C[G][t]$ by setting, identically on $G$, $q(t)=t^np(t^{-1})$.
A simple computation shows that, if we let $\chi$ be the character of $V$,
we have
  \[
  \frac d{dt}\ln q(t) = -\sum_{k\geq0}\psi^{k+1}(\chi)t^k.
  \]
It is clear that $p$ and $q$ have $1$ as a zero of the same multiplicity,
so that the function $d\in C[G]$ defined in \pref{the:theorem} is given by
  \[
  d=2n+2\res_1 \sum_{k\geq0}\psi^{k+1}(\chi)t^k,
  \]
where we have written $\res_1f$ the residue at $1$ of a function $f$
meromorphic in a neighborhood of~$1$.

This equation implies that the numbers $\dim_\C\HH^k(A_n{G})$
are determined by the character $\chi$ and the $\lambda$-ring structure on
$C[G]$. Since we are working over a field of characteristic zero, this last
structure is determined by the Adams operations, and these, in turn, depend
only on the power-maps, \ie~the maps induced on the set of the conjugacy
classes of $G$ by exponentiation. 

We thus see that if two groups $G$ and $G'$ are such that both their
character tables and power-maps coincide---such pairs are shown to exist in
\cite{Dade}---and we choose corresponding faithful representations,
which will of course have the same degree $n$, say, $\HH^*(A_n^G)$ and
$\HH^*(A_n^{G'})$ will be isomorphic as graded vector spaces. 

\paragraph In fact, these isomorphism can be taken to be multiplicative.
Let $\br G$ be the set of conjugacy classes of $G$, and, for $c\in\br G$,
let $\Hat c=\sum_{g\in c}g\in\C\,G$. The set $\{\Hat c\}_{c\in\br G}$ is a
basis of $\Z\,G$, the center of $\C\,G$.

%%%%%%%%%%%%%%%%%%%%%%%%%%%%%%%%%%%%%%%%%%%%%%%%%%%%%%%%%%%%%%%%%%%%%% 

\newpage

\paragraph For each $n\geq0$, let $S_n$ be the symmetric group on
$\{1,\dots,n\}$, and let $i_n:S_n\fn S_{n+1}$ be the standard injection,
under which $S_n$ fixes $n+1$;
let $S_\infty=\lim S_n$ be the injective limit, the
\emph{restricted symmetric group} on an countable infinite number of
letters. 

\paragraph A \emph{partition} $\l$ is a non-increasing sequence
of non-negative integers $(\l_i)_{i\geq 1}$ which eventually vanish; let
$\Pi$ be the set of all partitions. If $\l\in\Pi$, let $l(\l)$ stand for the number
of non-zero terms in $\l$, and define the \emph{weight} of $\l$ to be 
$|\l|=\sum_{i\geq1}\l_i$. Let $\Pi_n$ be the set of partitions of weight
$n$.

\paragraph If $\pi\in S_n$, the \emph{type} of $\pi$ is the partition $\rho(\pi)$ listing the
lengths of the cycles in a disjoint cycle decomposition of $\pi$; clearly,
$|\rho(\pi)|=n$. If $\rho(\pi)=(r_1,\dots,r_{l})$,
$\rho(i_n(\pi))=(r_1,\dots,r_{l},1)$, so that if we define the \emph{stable
type} of $\pi$ to be the partition $\rho^\#(\pi)=(r_1-1,\dots,r_{l}-1)$, we
see that this is compatible with the injections $i_n$, and in consequence
$\rho^\#$ is defined on $S_\infty$. For $\l\in\Pi$, we set $C_\l=\{\pi\in
S_\infty:\rho^\#(\pi)=\l\}$; it is easy to show that $\{C_\l\}_{\l\in\Pi}$
is precisely the decomposition of $S_\infty$ into conjugacy classes.

\paragraph For $n\geq0$, let $\Z_n=\Z(\ZZ\,S_n)$, and if $\l\in\Pi$,
$c_\l(n)=\sum_{g\in C_\l\cap S_n}g\in\Z_n$. Obviously we have $c_\l(n)=0$
if $|\l|+l(\l)>n$.

\paragraph\label{far} If $\l,\mu\in\Pi$, there are integeres $a_{\l\mu}^\nu(n)$,
$\nu\in\Pi$, such that
  \[
  c_\l(n)c_\mu(n)=\sum_{\nu\in\Pi}a_{\l\mu}^\nu(n) c_\nu(n),
  \]
with $a_{\l\mu}^\nu(n)=0$ if $|\nu|>|\l|+|\mu|$. These numbers have great
combinatorial interest, and explicit computations of 
specific cases can be found in the literature. In general, \cf~\cite{Farahat}, $a_{\l\mu}^\nu(n)$ depends
polynomially on $n$, and is actually independent of $n$ when
$|\nu|=|\l|+|\mu|$.

\paragraph\label{zinfty} Let $\B\subset\Q[t]$ of polynomials which take integer values on
integeres, and let $\Z_\infty$ be the (possibly non-associative)
$\B$-algebra which as a $\B$-module is free on the set
$\{c_\l\}_{\l\in\Pi}$, and whose product is determined by
  \[
  c_\l c_\mu=\sum_{\substack{\nu\in\Pi\\|\nu|\leq|\l|+|\mu|}}a_{\l\mu}^\nu\,c_\nu.
  \]
There are multiplicative $\ZZ$-linear morphisms $\Z_\infty\epi\Z_n$ given
by specialization $c_\l\mapsto c_\l(n)$ on the basis, and by evaluation of
polynomials at $n$ on the coefficients, which collect to give a
multiplicative map $\Z_\infty\fn\Pi_{n\geq0}\Z_n$, which turns out to be
injective. This implies that $\Z_\infty$ is an associative algebra.
Clearly, the kernel of $\Z_\infty\epi\Z_n$ is generated by those $c_\l$
such that $|\l|+l(\l)>n$ and polynomials in $\B$ which vanish on $n$.

\paragraph Let $F_n\Z_\infty$ be the $\B$-submodule spanned by the $c_\l$
with $|\l|\leq n$. We see at once that this defines an algebra filtration
$F_\*\Z_\infty$ on $\Z_\infty$, and, in view of the last statement of
\pref{far}, we can describe the associated graded algebra as follows: let
$\G$ be the $\ZZ$-algebra with $\ZZ$-basis $\{c_\l\}_{\l\in\Pi}$, and
product given by $c_\l c_\mu=\sum_{|\nu|=|\l|+|\mu|}a_{\l\mu}^\nu\,c_\nu$;
then $\gr\Z_\infty=\B\ot_\ZZ\G$.

\paragraph The epimorphisms $\Z_\infty\fn\Z_n$ of \pref{zinfty} are
compatible with the filtrations on the objects involved, so they give
epimorphisms on associated graded objects; but it is easy to see that these
are actually determined by epimorphisms $\G\epi\gr\Z_n$, given by 
$c_\l\mapsto c_\l(n)$.

\paragraph Let us write, for $\l=(r)$, $c_r=c_\l$. It is not difficult to
show that the set $\{c_r\}_{r\geq1}$ is algebraically independent in $\G$,
and generates it rationally.

\paragraph Let $\Lambda$ be the ring of symmetric functions with integer
coefficients on an countably infinite number of variables, and, for
each $i\geq0$, let $h_i$ be the $i$-th complete symmetric function, which
is the sum of all monomials of total degree $i$. It turns out that
$\Lambda=\ZZ[\{h_i\}_{i\geq 1}]$. If $\l=(r_1,\dots, r_l)\in\Pi$, we set
$h_\l=h_{r_1}\cdots h_{r_l}$.

\paragraph Define $u=\sum_{i\geq0}h_it^{i+1}\in\Lambda[[t]]$, and define elements
$h^*_i\in\Lambda$ so that $t=\sum_{i\geq0}h^*_iu^{i+1}$. Since the $h_i$
freely generate $\Lambda$, we can define a ring morphism
$\Psi:\Lambda\fn\Lambda$ with $\Psi(h_i)=h^*_i$. Now one can verify from
their definition that the $h^*_i$ are also algebraically independent and
generate $\Lambda$, so, in view of the symmetry of the contruction, we see
that $\Psi$ is actually an involution.

\paragraph Let $\br{{-}|{-}}$ be the bilinear form on $\Lambda$ with values in
$\ZZ$ such that $\br{h_\l|m_\mu}=\delta_{\l\mu}$, where the $m_\mu$ are the
monomial symmetric functions, obtained, for $\l=(r_1,\dots,r_l)$, by
symmetrization from $x_1^{r_1}\cdots x_l^{r_l}$. Define functions
$\{g_\l\}_{\l\in\Pi}$ such that $\br{g_\l,h^*_\mu}=\delta_{\l\mu}$.

\paragraph Now we can give a very concrete description of the ring $\G$
and, using the fact that the kernel of the epimorphisms
$\G\epi\gr\Z_n$ is easily identifiable, of the algebras $\HH(A_n^{S_n})$.
Indeed, it is a theorem proved in \cite{Mac}, example I.7.25, that the map
$\phi:\Lambda\fn\G$ such that $\phi(g_\l)=c_\l$ is a ring isomorphism.
This allows us to do explicit computations in the algebras
$\HH(A_n^{S_n})$, by translating the problem into one involving symmetric
polynomials. 

For each $n\geq0$, let $RS_n$ be the representation ring of $S_n$;
the sum $RS=\bigoplus_{n\geq0}RS_n$ is a strictly
commutative graded ring with product determined by its restrictions 
$RS_n\ot RS_m\fn RS_{n+m}$, given by $\chi\cdot\eta=\ind_{S_n\times
S_m}^{S_{n+m}}(\chi\times\eta)$. There is a very natural isomorphism of rings
$\Theta:RS\fn\Lambda$, essentially corresponding to taking the character of
representations. It would be interesting to be able to explicitly relate
elements and their products in $\HH^\*(A_n^{S_n})$ to actual
representations of the symmetric groups using the composition
$\phi\circ\Theta$ of isomorphisms at hand. Unfortunately, one cannot expect
to be able to restrict oneself to actual representations, since already
$g_{(1)}$ corresponds under $\Theta$ to $-1\in RS_1$, the opposite of
the trivial representation of $S_1$, which, of course, is only a virtual
representation. 

%%@@@%%%%%%%%%%%%%%%%%%%%%%%%%%%%%%%%%%%%%%%%%%%%%%%%%%%%%%%%%%%%%%%%% 

\end{document}